\documentclass[12pt,reqno]{article}
\usepackage[usenames]{color}
\usepackage{xcolor}
\usepackage{amssymb, amsmath, amsfonts, amsthm, amscd}
\allowdisplaybreaks
\usepackage{graphicx}
\usepackage[colorlinks=true,
linkcolor=webgreen,
filecolor=webbrown,
citecolor=webgreen]{hyperref}
\definecolor{webgreen}{rgb}{0,.5,0}
\definecolor{webbrown}{rgb}{.6,0,0}
\usepackage{fullpage}
\usepackage{float}
\usepackage{latexsym}
\usepackage{epsf}
\usepackage{mathrsfs}
\newcommand{\seqnum}[1]{\href{https://oeis.org/#1}{\rm \underline{#1}}}
\usepackage{tikz}
\usepackage{tkz-tab}
\usetikzlibrary{matrix,calc,fit}

\theoremstyle{plain}
\newtheorem{theorem}{Theorem}[section]

\theoremstyle{definition}

\theoremstyle{remark}

\begin{document}
\begin{center}
\vskip 1cm
{\LARGE\bf Normal Ordering and Bessel Numbers with Integral Operators} 
\vskip 1cm
\large
Abdelhay Benmoussa\\
School District Ansis, Unit Alanfoukt\\
Agadir Ida Outanane, Tiqqi 80324\\
Morocco\\
\href{mailto:abdelhay.benmoussa@taalim.ma}{\tt abdelhay.benmoussa@taalim.ma}
\end{center}
\vskip .2 in

\begin{abstract}
We derive a normal ordering formula for the operator \((xI)^n\), where \(I\) denotes the Volterra operator. The resulting coefficients are shown to coincide with the Bessel numbers. We also present two applications, along with a generalization of the main result.
\end{abstract}

\section{Introduction}

The normal ordering of differential operators is a classical topic in combinatorial analysis. In particular, the normal ordered form of the operator $(x\frac{d}{dx})^n$ is well known to involve coefficients given by the Stirling numbers of the second kind, a result that has appeared in various analytic and physical contexts \cite{MSS}.

In this paper, we investigate an alternative setting obtained by replacing the differentiation operator $\frac{d}{dx}$ with the Volterra operator
\[
(I f)(x) = \int_0^x f(t)\,dt.
\]
We show that the normal ordering of the operator $(xI)^n$ leads naturally to a new explicit formula whose coefficients are given by the Bessel numbers, defined by
\[
a(n,k)=[x^k]\,y_n(x)=\frac{(n+k)!}{2^k\,k!\,(n-k)!}
\]
where \(y_n(x)\) is the $n$-th Bessel polynomial \cite{Grosswald}.

Applying the formula to functions leads to closed-form expressions, which we illustrate through two examples.

Finally, we propose a generalization of the obtained formula, in which the coefficients are the Mansour--Shork numbers, providing a broader combinatorial framework for the integral-operator case.

\section{Main result}
\begin{theorem}\label{th}
Let \(n\ge1\). Then
\begin{equation}\label{000}
(xI)^n = \sum_{k=0}^{n-1} (-1)^k\, a(n-1,k)\, x^{\,n-k} I^{\,n+k}.
\end{equation}
\end{theorem}
\begin{proof}
Applying \((xI)^n\) to a function \(f\) gives 
\[
(xI)^n f(x) = x \int_0^x x_{n-1} \int_0^{x_{n-1}} \cdots x_1 \int_0^{x_1} f(t)\, dt \, dx_1 \cdots dx_{n-1}.
\]
This can be rewritten as
\[
(xI)^n f(x) = x \int_0^x f(t) 
\left( \iiint_{t \le x_1 \le \cdots \le x_{n-1} \le x} x_1 x_2 \cdots x_{n-1} \, dx_1 \cdots dx_{n-1} \right) dt.
\]
By induction, the inner integral evaluates to
\[
\iiint_{t \le x_1 \le \cdots \le x_{n-1} \le x} x_1 x_2 \cdots x_{n-1} \, dx_1 \cdots dx_{n-1} = \frac{(x^2 - t^2)^{n-1}}{2^{\,n-1}(n-1)!}.
\]
Thus we have
\[
(xI)^n f(x) = x \int_0^x \frac{(x^2 - t^2)^{\,n-1}}{2^{\,n-1}(n-1)!}\,f(t) \, dt \\
= \frac{x}{2^{\,n-1}(n-1)!} \int_0^x f(t) (x-t)^{\,n-1} (x+t)^{\,n-1} \, dt.    
\]
Expanding \((x+t)^{\,n-1}\) via the binomial theorem gives
\[
(x+t)^{\,n-1} = \sum_{k=0}^{n-1} \binom{n-1}{k} (t-x)^k (2x)^{\,n-1-k} = \sum_{k=0}^{n-1} (-1)^k \binom{n-1}{k} (2x)^{\,n-1-k} (x-t)^k.
\]
Substituting back and applying Cauchy's formula for repeated integration
\[
I^nf(x)=\frac{1}{(n-1)!}\int_0^x(x-t)^{n-1}f(t)\,dt,
\]
we obtain
\[
(xI)^n f(x) = \sum_{k=0}^{n-1} (-1)^k a(n-1,k)\, x^{\,n-k} I^{\,n+k} f(x),
\]
completing the proof.
\end{proof}
\section{Two applications}
\subsection{Power functions}
For \(\alpha>0\):
\begin{align*}
(xI)^{n+1}(t^{\alpha-1})(x)
&=x\int_0^x  \frac{(x^2 - t^2)^n}{2^n n!}t^{\alpha-1} \, dt \\
&= \frac{x}{2^n n!} \int_0^1 \big(x^2(1-u)\big)^n
   \big(x u^{1/2}\big)^{\alpha-1}\frac{x}{2}u^{-1/2}\,du
   \qquad (t=x\sqrt{u})\\
&= \frac{x}{2^n n!}\cdot \frac{x}{2}\, x^{2n} x^{\alpha-1}
   \int_0^1 (1-u)^n u^{(\alpha/2)-1}\,du \\
&= \frac{x^{\alpha+2n+1}}{2^{n+1}\, n!} \, B\Big(\frac{\alpha}{2},n+1\Big)\\
&= \frac{x^{\alpha+2n+1} \Gamma(\alpha/2)}{2^{n+1} \, \Gamma(n+\alpha/2+1)}\\
&=\frac{x^{\alpha+2n+1}}{\alpha(\alpha\,+2)\cdots(\alpha+2n)}.
\intertext{Series expansion gives}
(xI)^{n+1}(t^{\alpha-1})(x)
&=\sum_{k=0}^{n} (-1)^k a(n,k)\,
   x^{n+1-k}\,I^{n+1+k}(t^{\alpha-1})(x)\\
&= \sum_{k=0}^{n} (-1)^k a(n,k)\,
   x^{n+1-k} \int_0^x \frac{(x-t)^{n+k}}{(n+k)!} t^{\alpha-1}\,dt \\
&= \sum_{k=0}^{n} \frac{(-1)^k a(n,k)}{(n+k)!}\,
   x^{n+1-k} \int_0^1 x^{n+k}(1-u)^{n+k}
   x^{\alpha-1}u^{\alpha-1}\,x\,du
   \qquad (t=xu)\\
&= \sum_{k=0}^{n} \frac{(-1)^k a(n,k)}{(n+k)!}\,
   x^{2n+\alpha+1}
   \int_0^1 (1-u)^{n+k}u^{\alpha-1}\,du \\
&= \sum_{k=0}^{n} \frac{(-1)^k a(n,k)}{(n+k)!}\,
   x^{2n+\alpha+1} B(\alpha,n+k+1) \\
&= x^{2n+\alpha+1}\sum_{k=0}^{n}\,
   \frac{(-1)^k a(n,k)\Gamma(\alpha)}{\Gamma(\alpha+n+k+1)} \\
&= x^{\alpha+2n+1}\,
   \sum_{k=0}^{n}\frac{(-1)^k\,a(n,k)}{\alpha(\alpha+1)\cdots(\alpha+n+k)}
\end{align*}
Comparing the two expressions yields a generating function for the numbers $a(n,k)$:
\begin{equation}
\frac{1}{\alpha(\alpha\,+2)\cdots(\alpha+2n)}=\sum_{k=0}^{n}\frac{(-1)^k\,a(n,k)}{\alpha(\alpha+1)\cdots(\alpha+n+k)}
\end{equation}
\subsection{Exponential function}
\begin{align*}
(xI)^{n+1}(e^t)(x)
&=\sum_{k=0}^{n} (-1)^k a(n,k)\,
   x^{n+1-k}\,I^{n+1+k}(e^t)(x)\\
&= \sum_{k=0}^{n} (-1)^k a(n,k)\,
   x^{n+1-k} \int_0^x\frac{(x-t)^{n+k}}{(n+k)!}e^t\,dt.
\end{align*}
The integral evaluates as
\[
\int_0^x (x-t)^{n+k} e^t\,dt = e^x \, \gamma(n+k+1,x),
\]
where \(\gamma\) is the lower incomplete gamma function. Hence,
\begin{align*}
(xI)^{n+1}(e^t)(x)
&= e^x \sum_{k=0}^{n} (-1)^k a(n,k)\, x^{n+1-k} \,\frac{\gamma(n+k+1,x)}{(n+k)!}  \\
&= e^x \sum_{k=0}^{n} (-1)^k a(n,k)\, x^{n+1-k} 
   \left(1 - e^{-x}\sum_{j=0}^{n+k} \frac{x^j}{j!}\right) \\
&= e^x \sum_{k=0}^{n} (-1)^k a(n,k)\, x^{\,n+1-k}
   - \sum_{k=0}^{n}\sum_{j=0}^{n+k} (-1)^k a(n,k)\,\frac{x^{\,n+j+1-k}}{j!}\\
&= x^{\,n+1} e^{x}\,
y_n\!\left(-\frac1x\right)
-
\sum_{k=0}^{n}
(-1)^{\,n-k}\,
\frac{(2(n-k)-1)!!}{(2k)!!}\,
x^{\,2k+1}.  
\end{align*}
where the collapse of the double sum to the double factorial series was proved in \cite{RosengrenMO}. Evaluating at $x=1$ yields a Dobinski-type identity for the sequence
$a(n)=
y_n(-1)$:
\begin{equation}
a(n)
=
\frac1e\left(
 (xI)^{n+1}(e^t)(1)
+
\sum_{k=0}^{n}
(-1)^{\,n-k}\frac{(2(n-k)-1)!!}{(2k)!!}
\right).
\end{equation}
\section{Generalization}
Let $\mathscr{A}_{s;h}$ denote the algebra generated by two variables $U$ and $V$ satisfying the commutation relation
\[
UV - VU = h\,U^s,
\qquad h \in \mathbb{C} \setminus \{0\},\ s \in \mathbb{N}.
\]
Observe that 
\[
Ix - xI = -I^2,
\]  
which means the operators $x$ and $I$ generate $\mathscr{A}_{2;-1}$. Thus, Theorem \ref{th} is a special case of the following general normal ordering result for the word $(VU)^n$ in the algebra $\mathscr{A}_{s;h}$.

\begin{theorem}
For any integer $n \ge 1$, the normal ordered form of $(VU)^n$ in $\mathscr{A}_{s;h}$ is given by
\begin{equation}
(VU)^n
  = \sum_{k=1}^n \mathfrak{S}_{s;h}(n,k)\,
      V^{\,k}\, U^{\,s(n-k)+k},
\end{equation}
where $\mathfrak{S}_{s;h}(n,k)$ are the generalized Stirling numbers introduced by Mansour and Schork \cite{MSS}.    
\end{theorem}

\begin{proof}
We proceed by induction on $n$.

For $n=1$ we have
\[
(VU)^1 = VU.
\]
Since $\mathfrak{S}_{s;h}(1,1)=1$, the formula
\[
(VU)^n
= \sum_{k=1}^n \mathfrak{S}_{s;h}(n,k)\,
V^{k} U^{s(n-k)+k}
\]
holds for $n=1$.

Assume that for some $n \ge 1$,
\[
(VU)^n
= \sum_{k=1}^n \mathfrak{S}_{s;h}(n,k)\,
V^{k} U^{s(n-k)+k}.
\]

We compute
\begin{align*}
(VU)^{n+1}
&= (VU)^n (VU) \\[4pt]
&= \sum_{k=1}^n \mathfrak{S}_{s;h}(n,k)\,
V^{k} U^{s(n-k)+k} VU
\qquad
\text{(by the induction hypothesis)} \\[6pt]
&= \sum_{k=1}^n \mathfrak{S}_{s;h}(n,k)
\Bigl(
V^{k+1} U^{s(n-k)+k+1}
+ h\bigl(s(n-k)+k\bigr)
V^{k} U^{s(n-k)+k+s}
\Bigr) \\
&\hspace{7cm}
\text{(by Lemma 3.1 in \cite{MSS})} \\[6pt]
&= \sum_{k=2}^{n+1} \mathfrak{S}_{s;h}(n,k-1)\,
V^{k} U^{s(n+1-k)+k} + \sum_{k=1}^{n}
h\bigl(k+s(n-k)\bigr)\mathfrak{S}_{s;h}(n,k)\,
V^{k} U^{s(n+1-k)+k} \\[6pt]
&= \sum_{k=1}^{n+1}
\Bigl(
\mathfrak{S}_{s;h}(n,k-1)
+ h\bigl(k+s(n-k)\bigr)\mathfrak{S}_{s;h}(n,k)
\Bigr)
V^{k} U^{s(n+1-k)+k} \\[6pt]
&= \sum_{k=1}^{n+1}
\mathfrak{S}_{s;h}(n+1,k)\,
V^{k} U^{s(n+1-k)+k}
\qquad
\text{(by Equation (18) in \cite{MSS}).}
\end{align*}

This completes the induction and the proof.
\end{proof}

\bigskip
\hrule
\bigskip
\noindent 2020 \textit{Mathematics Subject Classification}: 11B83, 11C08, 26A36.

\noindent
\emph{Keywords:} Normal ordering, Bessel numbers, Iterated integrals, Operator calculus.
\bigskip
\hrule
\bigskip
\noindent (Concerned with sequences \seqnum{A000806}, \seqnum{A001497}, \seqnum{A001498}, \seqnum{A122848}, and \seqnum{A122850}.)
\bigskip
\hrule
\bigskip
\noindent
Return to \href{https://cs.uwaterloo.ca/journals/JIS/}{Journal of Integer Sequences home page}.
\vskip .1in

\begin{thebibliography}{99}
\bibitem{Grosswald} E. Grosswald, \textit{Bessel Polynomials}, Springer, 1978.
\bibitem{MSS}
T. Mansour, M. Schork and M. Shattuck, On a new family of generalized Stirling and Bell numbers, {\em Electron. J. Combin.} {\bf18} (2011) \#P77.
\bibitem{RosengrenMO}
H. Rosengren. 
\textit{Proving that a certain factorial double sum collapses to a double--factorial series}. 
MathOverflow posting (2025). Available at 
\url{https://mathoverflow.net/questions/499957}.
\end{thebibliography}
\end{document}